\documentclass[12pt]{huber_arxiv_article}

\usepackage[hmargin=1in,vmargin=1in]{geometry}
\usepackage{graphicx}
\usepackage{verbatim}
\usepackage{enumerate}
\usepackage{dsfont}
\usepackage{booktabs}
\usepackage{color}

\newcommand{\real}{\mathbb{R}}
\newcommand{\prob}{\mathbb{P}}
\newcommand{\mean}{\mathbb{E}}
\newcommand{\unifdist}{\textsf{Unif}}
\newcommand{\pois}{\textsf{Pois}}
\newcommand{\berndist}{\textsf{Bern}}
\newcommand{\ind}{{\mathds{1}}}
\newcommand{\tpa}{{\sffamily{TPA}}}

\usepackage[section]{algorithm}
\usepackage{algorithmic}

\newcommand{\DRAW}{\STATE {\bf draw }}
\newcommand{\uline}[1]{\underline{\makebox[12pt]{#1}}}

\newcommand{\linext}{\mathcal{L}}

\usepackage{amsmath,amsfonts,amsthm}

\theoremstyle{plain}
\newtheorem{theorem}{Theorem}

\newtheorem{lemma}{Lemma}

\theoremstyle{definition}

\newtheorem{definition}{Definition}

\begin{document}

\title{Using TPA to count linear extensions}

\maketitle

{\bfseries \sffamily Jacqueline Banks} \par
{\slshape University of California, Riverside} \par
{\slshape jbank003@student.ucr.edu}
\vskip 1em
{\bfseries \sffamily Scott M. Garrabrant} \par
{\slshape Pitzer College} \par
{\slshape scott@garrabrant.com} 
\vskip 1em
{\slshape }
{\bfseries \sffamily Mark L. Huber} \par
{\slshape Claremont McKenna College} \par
{\slshape mhuber@cmc.edu}
\vskip 1em
{\bfseries \sffamily Anne Perizzolo} \par
{\slshape Columbia University} \par
{\slshape aperizzolo11@students.claremont.edu}
{\slshape }

{\abstract A linear extension of a poset is a permutation of the elements
of the set that respects the partial order.  Let $\#(\linext)$ denote the number of
linear extensions.   It is a {\#}P complete problem
to determine $\#(\linext)$ exactly for an arbitrary poset, and so 
randomized approximation algorithms that draw randomly from the set of 
linear extensions are used.  In this work, the set of linear extensions is
embedded in a larger state space with a continuous parameter $\beta$.  The
introduction of a continuous parameter allows for the use of a more efficient
method for approximating $\#(\linext)$ called \tpa.
Our primary result is that it is possible to sample from this continuous
embedding in time that as fast or faster than the best known methods for 
sampling uniformly from linear extensions.    
For a poset containing $n$ elements, this
means we can approximate $\#(\linext)$ to 
within a factor of $1 + \epsilon$ with probability at least $1 - \delta$
using an expected number of random bits and comparisons in the poset 
which is at most 
$O(n^3 (\ln n)(\ln \#(\linext))^2 \epsilon^{-2} \ln \delta^{-1}).$}

\vspace{.1in}

\vspace{.1in}

\section{Introduction}

Consider the set $[n] = \{1,\ldots,n\}$.
A {\em partial order}
$\preceq$ is a binary relation on $[n]$ that is reflexive so $(\forall a \in [n])(a \preceq a)$, antisymmetric so 
$(\forall a,b \in [n])(a \preceq b \wedge
b \preceq a \Rightarrow a = b)$, and 
transitive so 
$(\forall a,b,c \in [n])(a \preceq b \wedge b \preceq c \Rightarrow a \preceq c)$.  A set equipped with a partial order is a partially ordered set, or 
{\em poset} for short.

The values $\{1,\ldots,n\}$ can be viewed as items.  
For a permutation $\sigma$, say that item $i$ has position $j$ if $\sigma(j) = i$.  Then a permutation 
respects the partial order $\preceq$ if whenever $i \preceq j$,
the position of item $i$ is less than the position of item $j$.  
Such a permutation is called a linear extension of the partial order.

\begin{definition}
  A permutation $\sigma$ is a {\em linear extension} of the
  partial order $\preceq$ if for all $i$ and $j$ in $[n]$, 
  $i \preceq j$ implies that $\sigma^{-1}(i) < \sigma^{-1}(j)$.
\end{definition}

For example, if the partial order states that $1 \preceq 3$ and 
$2 \preceq 4$, then the permutation 
$(\sigma(1),\sigma(2),\sigma(3),\sigma(4)) = (1,3,2,4)$ 
would be a linear extension.
However, $(4,1,3,2)$ would not since the position of 4 is 1, which is smaller
than the position of 2 which is 4.
This definition follows that of Karzanov and Khatchyan~\cite{karzanovk1991}. 
Note that some authors such as~\cite{alons2008} define the 
linear extension to be the permutation $\sigma^{-1}$ rather than 
$\sigma$.

Let $\linext$ denote the set of linear extensions of a particular poset.
Our goal here is to efficiently count the number of linear 
extensions, $\#(\linext)$.  Finding $\#(\linext)$ is a 
{\#}P complete problem~\cite{brightwellw1991} for general partial orders, 
and so instead of
an exact deterministic method, we develop a randomized 
approximation method.

There are many applications of this problem.  Morton 
et al.~\cite{mortonpssw2009} have shown that a particular type of 
convex rank test for nonparametric models can be reduced to counting
linear extensions.  Many data sets such as athletic competitions or
product comparisons do not have results for every possible pairing,
but instead have an incomplete set of comparisons.
Counting linear extensions can be used to develop estimates of the 
actual rank of the items involved (see~\cite{fishburng1975}.)

\paragraph{Previous results} 
Previous methods for this problem 
(\cite{karzanovk1991,huber2006b}) concentrated on 
sampling from the set of linear extensions where some 
of the permutation
values are fixed ahead of time.  Generating a single uniform
sample from the set of linear extensions takes 
$O(n^3 \ln n)$ expected number of random bits, using a number of 
comparisons that is at most the number of random bits~\cite{huber2006b}.
Using the self-reducibility method of Jerrum 
et al.~\cite{jerrumvv1986},
this can be used to estimate $\#(\linext)$ to within a factor of
$1 + \epsilon$ with probability at least $1 - \delta$ in time
$O(n^5 (\ln n)^3 \epsilon^{-2} \ln(1/\delta))$.

Here we take a different approach.  Instead
of sampling uniformly from the set of permutations, 
a weighted distribution
is used that has a parameter $\beta$.    
The weight assigned to an element varies continuously 
with $\beta$, and this
allows us to use a new method for turning samples from our weighted
distribution into an approximation for $\#(\linext)$ called the Tootsie Pop
Algorithm (\tpa).
The use of {\tpa} gives us an algorithm that is 
$O((\ln \#(\linext))^2 n^3 (\ln n) \epsilon^{-2} \ln(1/\delta))$.  In the
worse case, $\ln \#(\linext)$ is $O(n \ln n)$ and the complexity is the same
as the older algorithm, however, if $\#(\linext)$ is small compared to
$n!$, this algorithm can be much faster.  Even in the 
worst case, the constant hidden by the big-$O$ notation 
is much smaller for the new algorithm (see Theorem~\ref{THM:bits}
of Section~\ref{SEC:analysis}.)

\paragraph{Organization}
In the next section, we describe the self-reducibility method and {\tpa}
in detail.  Section~\ref{SEC:simple} illustrates the use of {\tpa} on
a simple example, and then Section~\ref{SEC:continuous} shows how it
can be used on the linear extensions problem by adding the appropriate
weighting.  Section~\ref{SEC:sampling} then shows how the non-Markovian
coupling from the past method introduced in~\cite{huber2006b} can also
be used for this new embedding, and Section~\ref{SEC:analysis} 
collects results concerning the running time of the procedure, including
an explicit bound on the expected 
number of random bits and comparisons used by the algorithm.

\section{The Tootsie Pop Algorithm}
\label{SEC:tpa}

In~\cite{jerrumvv1986}, Jerrum et al.~noted that for 
self-reducible problems, an algorithm for generating from
a set could be used to 
build an approximation algorithm for finding the size of the set.
Informally, a problem is self-reducible if the set of solutions
can be partitioned into the solutions of smaller instances of the 
problem (for precise details, see~\cite{jerrumvv1986}.)

For example, in linear extensions once the value of 
$\sigma(n)$ is determined,
the problem of 
drawing $\sigma(1),\ldots,\sigma(n-1)$ is just a smaller linear
extension generation problem.  

While a theoretical tour de force, as a practical matter using
self-reducibility to build algorithms is difficult.  The output of 
a self-reducibility 
algorithm is a scaled product of binomials, not the easiest
distribution to work with or analyze precisely.

The Tootsie Pop Algorithm (\tpa)~\cite{huber2010b} 
is one way to solve this difficulty.
Roughly speaking, {\tpa} begins with a large set (the {\em shell}) 
containing a smaller set (the {\em center}).  At each step, 
{\tpa} draws a sample $X$ randomly 
from the shell, and reduces the shell as much as possible while
still containing $X$.  The process then repeats, drawing
samples and contracting the shell.  This continues until the sample
drawn lands in the center.
The number of samples drawn before one falls in the center has a 
Poisson distribution, with parameter equal to the natural logarithm
of the ratio of the size of the shell to the center.

To be precise, {\tpa} requires the following ingredients 
\begin{enumerate}[(a)]
\item{A measure space $(\Omega,{\cal F},\mu)$.}
\item{Two finite 
measurable sets $B$ and $B'$ satisfying $B' \subset B$.  The
set $B'$ is the {\em center} and $B$ is the {\em shell}.}
\label{ing:sets}
\item{A family of nested sets $\{A(\beta):\beta \in \real\}$ 
      such that $\beta < \beta'$ implies 
$A(\beta) \subseteq A(\beta')$.  Also
$\mu(A(\beta))$ must be a continuous
function of $\beta$, and 
$\lim_{\beta \rightarrow -\infty} \mu(A(\beta)) = 0.$}
\item{Special values $\beta_B$ and $\beta_{B'}$ that satisfy
$A(\beta_B) = B$ and $A(\beta_{B'}) = B'$.}
\end{enumerate}

With these ingredients, {\tpa} can be run as follows.  
\algsetup{indent=2em}
\begin{algorithm}[ht]
\caption{\quad TPA$(r,\beta_B,\beta_{B'})$}\label{ALG:TPA}
\begin{algorithmic}[1]
\REQUIRE Number of runs $r$, initial index $\beta_B$, final index
          $\beta_{B'}$ 
\ENSURE  $\hat L$ (estimate of $\mu(B)/\mu(B')$)
\STATE $k \leftarrow 0$
\FOR{$i$ from $1$ to $r$}
  \STATE $\beta \leftarrow \beta_{B},$ $k \leftarrow k - 1$
  \WHILE{$\beta > \beta_{B'}$}
    \STATE $k \leftarrow k + 1$, $X \leftarrow \mu(A(\beta))$,
           $\beta \leftarrow \inf\{\beta' \in [\beta_{B'},\beta_B]:
              X \in A(\beta')\}$ 
  \ENDWHILE
\ENDFOR
\STATE $\hat L \leftarrow \exp(k / r)$
\end{algorithmic}
\end{algorithm}

Let $A = \ln(\mu(B)/\mu(B'))$, so that $\exp(A)$ is what we
are trying to estimate.  Then
each run through the for loop in the algorithm requires on average
$A + 1$ samples, making the total expected
number of samples $r (A + 1)$.  The value of $k$ in line 7 of the algorithm
is Poisson distributed with parameter $r A$.
This means that $r$ should be set to about $A$ so that 
$k / r$ is tightly concentrated around $A$.

But we do not know $A$ ahead of time!  This leads to the need for a 
two-phase algorithm.  In the first phase $r$ is set to be large enough
to get a rough approximation of $A$, and then in the second phase
$r$ is set based on our estimate from the first run.  That is:
\begin{enumerate}
\item{Call {\tpa} with $r_1 = 2 \ln(2/\delta)$ to obtain $\hat L_1$, and
set $\hat A_1 = \ln(\hat L_1)$.}
\item{Call {\tpa} with $r_2 = 2 (\hat A_1 + \sqrt{\hat A_1} + 2) 
 [\ln(1 + \epsilon)^2 - \ln(1 + \epsilon)^3]^{-1} 
  \ln(4/\delta)$ to obtain the final
estimate.}
\end{enumerate}
The result is output $\hat L_2$ that is within a factor of 
$1 + \epsilon$ of $\#(\linext)$ with probability at least $1 - \delta$.
This is shown in Section~\ref{SEC:analysis}.

\section{Continuous embedding:  simple example}
\label{SEC:simple}

To illustrate {\tpa} versus the basic self-reducibility approach,
consider a simple problem that will serve as a building block for 
our algorithm on linear extensions later.  In this problem, 
we estimate the size of the 
set $\{1,2,\ldots,n\}$ given the ability to draw samples
uniformly from $\{1,2,\ldots,b\}$ for any $b$.  

In the self-reducibility approach, 
begin by setting $\beta_1 = \lceil n / 2 \rceil$ and drawing
samples from $\{1,\ldots,n\}$.  Count how many fall into
$\{1,\ldots,\beta_1\}$ and use this number $\hat a_1$ 
(divided by the number
of samples) as an estimate of $\beta_1/n$.  Now repeat, letting
$\beta_2 = \lceil \beta_1 / 2 \rceil$ and estimating 
$\hat a_2 = \beta_2 / \beta_1$ until $\beta_k = 1$.  Note that
\[
\mean[\hat a_1 \hat a_2 \cdots \hat a_{k-1}]
 = \frac{\beta_1}{n} \frac{\beta_2}{\beta_1} \cdots 
   \frac{\beta_{k}}{\beta_{k-1}} = \frac{\beta_k}{n}.
\]
Since the final estimate $\hat a$ of $\frac{1}{n}$
is the product of $k - 1$ estimates, Fishman called this algorithm
the {\em product estimator}~\cite{fishman1994}.  The problem with
analyzing the output of the product estimator, is that it is the
product of $k$ scaled binomials.  

To use {\tpa} on this problem, it needs to be embedded in a continuous
setting.  Consider the state space $[0,n]$.  The family of sets
needed for {\tpa} will be $[0,\beta]$, where $\beta_B = n$ and 
$\beta_{B'} = 1$.  This makes the ratio of the measure of $[0,\beta_B]$
to $[0,\beta_{B'}]$ equal to $n$.

Note that you can draw uniformly from $[0,\beta]$ in the following
two step fashion.  First draw 
$X \in \{1,2,\ldots,\lceil \beta \rceil\}$ so that 
$\prob(X = i) = 1/\beta$ for $i < \beta$ and 
$\prob(X = \beta) = 
(1 + \beta - \lceil \beta \rceil) / \beta$.
If $X < \beta$, draw $Y$ uniform on $[0,1]$, otherwise draw 
$Y$ uniform on $[0,1 + \beta - \lceil \beta \rceil]$.  
The final draw is $W = X - 1 + Y$.

{\tpa} starts with $\beta_0 = n$, then draws $W$ as above.
The infimum over all $\beta$ such that $W \in [0,\beta]$ is just
$\beta = W$.  So $\beta_1$ just equals $W$.  
Next, redraw $W$ from $[0,\beta_1]$.  Again, the infimum of 
$\beta$ satisfying $W \in [0,\beta]$ is just $W$, so $\beta_2$
equals this new value of $W$.  

This process repeats until 
$W$ falls into $[0,1]$.  The estimate $k$ for $\ln n$ is just
the number of steps needed before the final step into $[0,1]$.
Note that $k$ can equal 0 if the very first step lands in $[0,1]$.
This random variable $k$ will be Poisson distributed with parameter 
$\ln n$.  Recall that the sum of Poisson random variables is also
Poisson with parameter equal to the sum of the individual parameters, 
so repeating the process $r$ times and summing the results
yields a Poisson random variable with parameter $r \ln n$.  Dividing
by $r$ and exponentiating then yields an estimate of $n$.

\section{Continuous embedding:  linear extensions}
\label{SEC:continuous}

This approach can be extended to the problem of linear extensions 
by adding an auxilliary random variable.  First we define a distance
between an arbitrary permutation and a home linear extension.

Note that we can assume without loss of generality that
$(1,2,\ldots,n)$ is a valid linear extension, otherwise, simply relabel
the items so that it is.  Then say that $i$ is the {\em home position} 
of item $i$.

Let $a^+ = \max\{a,0\}$.  In linear extension $\sigma$, item $j$ has
position $\sigma^{-1}(j)$.
Define the distance
from item $j$ to its home position to be 
$(\sigma^{-1}(j) - j)^+$.  The maximum of
these distances over all items is the distance from $\sigma$ to the home position
That is, let
\[
d(\sigma,(1,2,\ldots,n)) = \max_j 
  (\sigma^{-1}(j) - j)^+.
\]

If the distance is 0, then no element $i$ is to the right of 
the home position.  The only way that can happen is if 
$\sigma(i) = i$ for all $i$.

Right now, the distance is discrete, falling into $\{0,1,2,\ldots,n-1\}$,
and all linear extensions are equally likely.
To finish the continuous embedding, it is necessary to change from
a uniform distribution to one where some linear extensions are
more likely than others.  

Let $\beta \in [0,n-1]$.
Suppose that item $i$ is farther than $\beta$ to the right
of its home position.  Such an item has weight 0.  
If item $i$ is closer than $\lceil \beta \rceil$ to its home position,
it has weight $1$.  If item $i$ is exactly distance $\lceil \beta \rceil$
from its home position, then it receives weight equal to 
$1 + \beta - \lceil \beta \rceil$.

For instance, for $\sigma = (1,4,3,2)$,
and $\beta = 3$,
the weight of items $1$, $4$, and $3$, is 1, while the weight of item 
$2$ is $1 + 3 - 3 = 1$.
If $\beta$ falls to $2.3$ then the weight of item 
$2$ drops to $1 + 2.3 - 3 = 0.3$.

Let the weight of a linear extension be the product of the weights
of each of its items.  That is,
\begin{equation}
\label{EQN:weights}
w(\sigma,\beta) = \prod_{i \in [n]} w_i(\sigma,\beta), \text{ where }
  w_i(\sigma,\beta) = ((1 + \beta - \lceil \beta \rceil) 
   \ind(\sigma^{-1}(i) - i = \lceil \beta \rceil)
   + \ind(\sigma^{-1}(i) - i < \lceil \beta \rceil)).
\end{equation}

In other words, when $\beta$ is an integer, all weights are either 0 or 1,
and are 1 if and only if all items are at distance at most $\beta$ to the
right of their home position.  When
$\beta$ is not an integer, then all items must be at most 
$\lceil \beta \rceil$ distance from home, and every item whose distance from home
equals $\lceil \beta \rceil$ receives a penalty factor equal to the 
fractional part of $\beta$.

Note that $w(\sigma,\beta)$ is an increasing function of $\beta$.

Suppose $X$ is a random element of $\linext$ where 
$\prob(X = \sigma) \propto w(\sigma,\beta)$.  Let $\unifdist(\Omega)$
denote the uniform distributon over the set $\Omega$.  Given $X$, create the 
auxiliary variable $Y$ as
$[Y|X] \sim \unifdist([0,w(\sigma,\beta)])$.  Let 
$A(\beta) = \{(x,y):x \in \linext,y \in [0,w(x,\beta)]\}$.

It is upon these sets $A(\beta)$ that \tpa{} can be used.
Here $A(n-1) = \linext \times [0,1]$ is the shell and
$A(0) = \{(1,2,\ldots,n)\} \times [0,1]$ is the center.
Then since $w(\sigma,\beta)$ is an increasing function of $\beta$,
for $\beta' \leq \beta$, $A(\beta') \subseteq A(\beta)$.

\section{Sampling from the continuous embedding}
\label{SEC:sampling}

For the continuous embedding 
to be useful for \tpa, it must be possible to sample from the 
set of linear extensions 
with the weight function given in~\eqref{EQN:weights}.  Once the linear
extension $X$ has been created, sampling the $Y$ to go along with it
is simple.

To sample from the set of weighted linear extensions, first consider
a Markov chain whose stationary distribution matches the target
distribution.  This is done by using a Metroplis-Hastings approach.
The proposal chain works as follows.  
With probability 1/2, the chain just stays where it is.  With probability
$1/2$, a position $i$ is chosen uniformly from
$\{1,\ldots,n-1\}$.  If such a transposition obeys the partial
order and does not move an item more than $\lceil \beta \rceil$ 
to the right of its home position, it is the proposed move.

If the proposal is to transpose the items, then one item might have
acquired a weight factor of $1 + \beta - \lceil \beta \rceil$ if
it moves to be exactly $\lceil \beta \rceil$ from its home position.
So we only accept such a move with probability 
$1 + \beta - \lceil \beta \rceil$.

This is encoded in Algorithm~\ref{ALG:chain}.
\algsetup{indent=2em}
\begin{algorithm}[h!]
\caption{\quad {\tt ChainStep}($\sigma, i, C_1,C_2$)}\label{ALG:chain}
\begin{algorithmic}[1]
\REQUIRE current linear extension Markov chain state $\sigma$
\ENSURE  next linear extension Markov chain state $\sigma$
\STATE{$d \leftarrow i + 1 - \sigma(i)$}
\IF{$C_1 = 1$ and not $\sigma(i) \preceq \sigma(i + 1)$ and 
 $d \leq \lceil \beta \rceil$}
  \IF{$d < \beta$
      or $C_2 = 1$}
    \STATE{$a \leftarrow \sigma(i + 1)$, $\sigma(i + 1) \leftarrow \sigma(i)$,
           $\sigma(i) \leftarrow a$}
  \ENDIF
\ENDIF
\end{algorithmic}
\end{algorithm}

Write $B \sim \berndist(p)$ if $\prob(B = 1) = p$ and $\prob(B = 0) = 1-p$.
Then with the appropriate choice of random inputs, Algorithm~\ref{ALG:chain}
has the distribution on linear extensions with probabilities proportional
to $w(\cdot,\beta)$ as a stationary distribution.

\begin{lemma}
  For $\sigma \sim w(\sigma,\beta)$, $i \sim \unifdist(\{1,\ldots,n-1\})$,
  $C_1 \sim \berndist(1/2)$, 
  $C_2 \sim \berndist(1 + \beta - \lceil \beta \rceil)$.  Then
  {\tt ChainStep$(\sigma,i,C_1,C_2) \sim w(\sigma,\beta)$}.
\end{lemma}

\begin{proof}
  This follows from the reversibility (see for instance~\cite{resnick1992})
  of the Markov chain with respect to $w$.
\end{proof}

From this chain, it is possible to build a method for obtaining
samples exactly from the target distribution.  The method of 
coupling from the past (CFTP) was developed by Propp and 
Wilson~\cite{proppw1996} to draw samples exactly from the stationary
distribution of Markov chains.  For this problem, an extension called
non-Markovian CFTP~\cite{huber2006b} is needed.

The method works as follows.  First, a bounding chain~\cite{huber2004a}
is constructed for the chain in question.  A bounding chain is an
auxiliary chain on the set of subsets of the original state space.
That is, $\Omega_{\textrm{bound}} = 2^{\Omega}$, where $\Omega$ is the 
state space of the original chain.  Moreover, 
there is a coupling between
the original chain $\{\sigma_t\}$ and 
the bounding chain $\{S_t\}$ such 
that $\sigma_t$ evolves according to the kernel of the original bounding
chain, and 
$\sigma_t \in S_t \rightarrow \sigma_{t+1} \in S_{t+1}$.  

For us, the state of the bounding chain is indexed by a 
vector $B \in \{1,\ldots,n, \theta\}^n.$  Let
\[
S(B) = \{\sigma:(\forall i)((B(j) = i) \wedge (\sigma(j') = i) \Rightarrow 
    j' \leq j)\}
\]
For instance, if $B(3) = 4$, then $\sigma \in S(B)$ requires that 
$\sigma(1) = 4$ or $\sigma(2) = 4$ or $\sigma(3) = 4$.  In this setup
$\theta$ is a special symbol:  if $B(i) = \theta$, then there is no
restriction on $\sigma$ whatsoever.  To visualize what is happening with
the state and bounding state, it will be useful to have a pictorial 
representation.  For instance, if $\sigma = (4,2,3,1)$ and 
$B = (\theta,4,3,\theta)$ this can be represented by:
\[
\uline{4}|_\theta\uline{2}|_4\uline{3}|_3 \uline{1}|_\theta.
\]
The bounding state works by keeping track of the right most position 
of the item in the underlying
state.  If $B(i) = a$, say that bar $|_a$ is at position $i$. 
To be a bounding state, if bar $|_a$ is at position $i$, then item
$a$ must be at a position in $\{1,2,\ldots,i\}$.

Now suppose there is a single $|_1$ at the rightmost position and
all other positions contain $|_\theta$.  Then this state
$B = (\theta,\ldots,\theta,1)$ bounds all permutations.

Next, suppose that there are no $|_\theta$ anywhere in the bounding state.
For instance $B = (2,4,1,3)$.  Let $x$ be a state bounded by $B$.  
Then $B(1) = 2$ means that item 2 in in position 1.  $B(2) = 4$ 
means that item 4 is in position 1 or 2.  But item 2 is in position 1,
so 4 must be in position 2.  Similarly, item 1 must be in position 3 and
item 3 must be in position 4.  In other words, if no component of $B$
is labeled $\theta$, then $S(B) = \{B\}$.  In our example
  \[
  \uline{2}|_2\uline{4}|_4\uline{1}|_1 \uline{3}|_3.
  \]


We are now ready to state the procedure for updating the current state
and the bounding state simultaneously.  This operates as in
Algorithm~\ref{ALG:bounding}.
Note that if the inputs to the Algorithm have
$i \sim \unifdist(\{1,2,\ldots,n\})$ and $C_1 \sim \berndist(1/2)$, then 
the state $\sigma$ is updated using the same probabilities as the 
previous chain step.
The key difference between how $\sigma$ and $B$ are updated 
is that if $\sigma(i) = B(i + 1)$, then 
$B$ is updated using $C_3 = 1 - C_1$, otherwise $C_3 = C_1$.
In any case, since $C_1 \sim \berndist(1/2)$, 
$C_3 \sim \berndist(1/2)$ as well.

\algsetup{indent=2em}
\begin{algorithm}[h!]
\caption{\quad {\tt BoundingChainStep}($\sigma, B, i, C_1, C_2$)}\label{ALG:bounding}
\begin{algorithmic}[1]
\REQUIRE current state and bounding state $(\sigma,B)$
\ENSURE  next state and bounding state $(\sigma,B)$
\STATE{$C_3 \leftarrow (1 - C_1) \ind(\sigma(i) = B(i+1))
 + C_1 \ind(\sigma(i) \neq B(i+1))$}
\STATE{$\sigma \leftarrow {\tt ChainStep}(\sigma,i,C_3,C_2)$}
\STATE{$B \leftarrow {\tt ChainStep}(B,i,C_1,C_2)$}
\IF{$B(n) = \theta$}
  \STATE{$B(n) \leftarrow 1 + \#\{j:B(j) \neq \theta\}$}
\ENDIF
\STATE{Return $(\sigma,B)$}
\end{algorithmic}
\end{algorithm}

Note that $\sigma$ is being updated as in Algorithm~\ref{ALG:chain}.
The only different is the bounding state update.  First, note that
if $i + \lceil \beta \rceil \leq n$, then the rightmost position that
item $i$ can be is $i + \lceil \beta \rceil$.  Hence there should be 
a $|_i$ at position $i + \lceil \beta \rceil$ or less.

\begin{definition}
  A bounding state $B$ is {\em $\beta$-tight} if for all items $i$ 
  with $i + \lceil \beta \rceil \leq n$, there exists 
  $j \leq i + \lceil \beta \rceil$ such that $B(j) = i$.
\end{definition}
Our main result is:
\begin{theorem}
If $\sigma \in S(B)$ for $B$ a $\beta$-tight bounding state, 
then running one step of Algorithm~\ref{ALG:bounding}
leaves $\sigma \in S(B)$ regardless of the inputs $i$, $C_1$ and $C_2$.
\end{theorem}

\begin{proof}
When $C_1 = C_3 = 0$, neither the $\sigma$ state or the 
$B$ state changes, and so the result is trivially true.

Write $(\phi(\sigma),\phi(B))$ for the output of the algorithm, supressing
the dependence on $i$, $C_1$, and $C_2$.  Given permutation $x$, write
$t(x,i)$ for the permutation where $x(i)$ and $x(i+1)$ have been
transposed.

In order for $\phi(\sigma) \notin \phi(B)$, there must be an item 
$a$ that moves to the right of the bar $|_a$.  If there is no $|_a$
in the bounding state (so there does not exist $j$ with $B(j) = a$)
then this trivally cannot happen.

Both bars and items can each move at most one step to the right or left.
So if either the position of 
$a$ is two or more to the left of the position of the bar $|_a$, or
there is no bar $|_a$ in the bounding state, then this also cannot
happen.

With that in mind, suppose $a$ is exactly 
one position to the left of the bar $|_a$.  Then
the only way that $a$ and $|_a$ could cross is if 
$\sigma^{-1}(a) = i$, $B^{-1}(a) = i+1$, and
$\phi(\sigma) = t(\sigma,i)$ and $\phi(B) = t(\sigma,i)$.  
But when $\sigma(i) = a = B(i+1)$, $C_3 = 1 - C_1$, so either
$\phi(\sigma) = \sigma$ or $\phi(B) = B$.  So this bad case cannot occur.

Suppose $a$ and $|_a$ are at the same position.  If that position is 
$i$, then since $\phi(\sigma)$ and $\phi(B)$ are using the same inputs 
and the weight factor incurred by moving $a$ to position $i+1$ is the 
same for both, 
either both use the transpose or neither do.  So either way no violation
occurs.

The final possibility to consider is that $\sigma(i+1) = B(i+1) = a$.
Is it possible for $|_a$ to move one position to the left while
$a$ stays at position $i + 1$?  Fortunately, the answer is once again 
no.  If $C_1 = 0$, then $C_3 = 0$, so both $\phi(\sigma) = \sigma$
and $\phi(B) = B$, so there is nothing to show.

Suppose $C_1 = C_3 = 1$.  Now consider the value of 
$\sigma(i)$.  If $i = \sigma(i) + \lceil \beta \rceil$, then 
the transpose operation on $\sigma$ would move item $\sigma(i)$ too
far to the right, and so $\phi(\sigma) = \sigma$.  But in this
case, since $B$ is $\beta$-tight,
$i = B(i) + \lceil \beta \rceil$ as well, and so $\phi(B) = B$.

Similarly, if $i = \sigma(i) + \lceil \beta \rceil - 1$, then 
either $B(i) = \sigma(i)$ or $B(i+1) = \sigma(i)$.  But the 
$B(i+1) = \sigma(i)$ case was dealt with earlier, leaving again
that $B(i) = \sigma(i)$.  So now if $C_2 = 1$ then
$\phi(\sigma) = t(\sigma,i)$ and $\phi(B) = t(B,i)$, and if 
$C_2 = 0$ then $\phi(\sigma) = \sigma$ and $\phi(B) = B$.  Either
way, they both move together.

If $i < \sigma(i) + \lceil \beta \rceil - 1$, then 
$\phi(\sigma) = t(\sigma,i)$, so there can never be a violation.

Hence in all cases $\sigma \in S(B) \Rightarrow \phi(\sigma) \in S(\phi(B))$.

\end{proof}

So if $\sigma$ is bounded by a $\beta$-tight $B$, 
it will still be bounded after taking
one step in the bounding chain step.  With this established, samples
from the target distribution can be generated as 
in Algorithm~\ref{ALG:generate}~\cite{huber2006b,proppw1996}
using non-Markovian CFTP.

\algsetup{indent=2em}
\begin{algorithm}[h!]
\caption{\quad {\tt Generate}($t$)}\label{ALG:generate}
\begin{algorithmic}[1]
\REQUIRE $t$ number of steps to use to generate a sample
\ENSURE  $\sigma$ drawn from the weighted distribution
\STATE{$\sigma \leftarrow (1,2,\ldots,n)$, $B \leftarrow
  (\theta,\ldots,\theta)$}
\FOR{$i$ from $1$ to $n - \lceil \beta \rceil$}
  \STATE{$B(i + \lceil \beta \rceil) \leftarrow i$}
\ENDFOR
\STATE{$B_0 \leftarrow B$}
\FOR{$j$ from $1$ to $t$}
  \DRAW{$i(j) \leftarrow \unifdist([n-1])$, 
      $C_1(j) \leftarrow \berndist(1/2)$,
      $C_2(j) \leftarrow \berndist(1 + \beta - \lceil \beta \rceil)$}
  \STATE{$(\sigma,B) \leftarrow 
     {\tt BoundingChainStep}(\sigma,B,i(j),C_1(j),C_2(j))$}
\ENDFOR
\IF{for all $i$, $B(i) \neq \theta$}
  \STATE{$\sigma \leftarrow B$}
\ELSE
  \STATE{$\sigma \leftarrow {\tt Generate}(2t)$, $B \leftarrow B_0$}
  \FOR{$j$ from $1$ to $t$}
  \STATE{$(\sigma,B) \leftarrow 
     {\tt BoundingChainStep}(\sigma,B,i(j),C_1(j),C_2(j))$}
  \ENDFOR
\ENDIF
\end{algorithmic}
\end{algorithm}

\section{TPA for linear extensions}

Now \tpa{} can be applied to linear extensions.  In the presentation
earlier, given $X \sim w(\cdot)$, 
a single random variable $[X|Y] \sim \unifdist([0,w(X)])$
was used to make the joint distribution uniform.  Since $w(x)$
has a product form, however, it makes things easier to generate
$n$ different auxiliary random variables $Y_1,\ldots,Y_n$ to make
it work.  If $w(x) = \prod_{i \in [n]} w_i(x)$, let
each $[Y_i|X]$ be independent and $\unifdist([0,w_i(X)]$.

Suppose that $Y_2 = 0.3$.  Then if item $2$ is 3 units to the right
of its home position, then that implies that $\beta \geq 2.3$.  If
item $2$ is 2 units to the right of its home position then 
$\beta \geq 1.3$, if it is 1 unit to the right of home then
$\beta \geq 0.3$.  Finally, if item $2$ is at home then 
$\beta \geq 0$.  In general, for item $j$ exactly
$X^{-1}(j) - j$ to the right of its home position, 
$\beta \geq b_i = (X^{-1}(j) - j  - 1 + Y_j)\ind(X^{-1}(j) - j > 0)$.

So that means that the next value of $\beta$ should be 
equal to the largest value of $b_i$.  Since this minimum is
taken over all items $i$, and $X$ is a permutation, the new value
of $\beta$ can be set to be
\[
\max_j (j - X(j) - 1 + Y_j)\ind(j - X(j) > 0).
\]

\algsetup{indent=2em}
\begin{algorithm}[ht]
\caption{\quad {\tt TPALinearExtensions}$(r)$}\label{ALG:TPALinearExtensions}
\begin{algorithmic}[1]
\REQUIRE Number of runs $r$
\ENSURE  $\hat L$ (estimate of $\#(\linext)$)
\STATE $k \leftarrow 0$
\FOR{$i$ from $1$ to $r$}
  \STATE$ \beta \leftarrow n - 1,$ $k \leftarrow k - 1$
  \WHILE{$\beta > 0$}
    \STATE $k \leftarrow k + 1$
    \STATE $X \leftarrow {\tt Generate}(1)$
    \FOR{$j$ from 1 to $n$}
       \STATE{{\bf draw} $Y_j \leftarrow \unifdist([0,w_j(X,\beta)])$}
       \STATE{{\bf let} $\beta_j \leftarrow 
          (X^{-1}(j) - j - 1 + Y_j)\ind(X^{-1}(j) - j > 0)$}
    \ENDFOR
    \STATE $\beta \leftarrow \max_j \beta_j$
  \ENDWHILE
\ENDFOR
\STATE $\hat L \leftarrow \exp(k / r)$
\end{algorithmic}
\end{algorithm}

\section{Analysis}
\label{SEC:analysis}

In this section we prove several results concerning the running time
of the procedure outlined in the previous section.

\begin{theorem}
\label{THM:nmcftp}
The non-Markovian coupling from the past in 
Algorithm~\ref{ALG:generate}
requires an expected number of random bits bounded by 
$4.3n^3(\ln n)(\lceil \log_2 n \rceil + 3)$ 
and a number of comparisons bounded by
$8.6 n^3 \ln n$.
\end{theorem}

\begin{theorem}
\label{THM:tpa}
For $\epsilon \leq 1$, the two-phase {\tpa} approach outlined at the end of 
Section~\ref{SEC:tpa} generates output $\hat L_2$ such that
\[
\prob((1 + \epsilon)^{-1} \leq \hat L_2 / \#(\linext) \leq 
 1 + \epsilon) \geq 1 - \delta.
\]
\end{theorem}

\begin{theorem}
\label{THM:bits}
The expected number of random bits needed to approximate 
$\#(\linext)$ to within a factor of $1 + \epsilon$ with probability at least
$1 - \delta$ is bounded above by 
\[
4.3n^3(\ln n)(\lceil \log_2 n\rceil + 3)
 [2 (A + 1) \ln (2/\delta) + (A + 1) (A + \sqrt{A} + 2)
  (\ln(1 + \epsilon)^2 - \ln(1 + \epsilon)^3) 
  \ln(4/\delta)].
\]

\end{theorem}

\begin{proof}[Proof of Theorem~\ref{THM:nmcftp}]
Lemma 10 of~\cite{huber2006b} showed that when there is no $\beta$
parameter, the expected number of steps taken by non-Markovian CFTP
was bounded above by $4.3n^3 \ln n.$

So the question is:  once the $\beta$ parameter falls below $n$,
does the bound still hold?  The bound was derived by considering
how long it takes for the $|_\theta$ values in the bounding state
 to disappear.  Each time
a $|_\theta$ reaches position $n$, it is removed and replaced by 
something of the form $|_a$.  When all the $|_\theta$ disappear,
the process in Algorithm~\ref{ALG:generate} terminates.

When there is no $\beta$, the probabilities that a particular 
$|_\theta$ bound moves to the left or the right are equal:  
both $1/(2n)$.  (This does not apply when
the bound is at position 1, in which case the bound cannot move
to the left.)  The result in~\cite{huber2006b} is really a bound on 
the number of steps in a simple random walk necessary for the $|_\theta$
bounds to all reach state $n$.

Now suppose that $\beta \in (0,n)$.  The probability that a $|_\theta$
bound moves to the right is still $1 /(2n)$, but now consider when the 
state is of the form $\ldots \uline{}|_a\uline{}|_\theta \ldots$. 
For $|_\theta$ to move left the $|_a$ has to move right, and this could
occur with probability $(1 + \beta - \lceil \beta \rceil)/(2n)$.
That is, with $\beta \in (0,n)$, the chance that the $|_\theta$ 
moves left can be below $1/(2n)$.

This can only reduce the number of moves necessary for the $|_\theta$
bounds to reach the right hand side!  That is, the random variable that
is the number of steps needed for all the $|_\theta$ bounds to reach position
$n$ and disappear is dominated by the same random variable for 
$\beta = n$.  Hence the bound obtained by
Lemma 10 of~\cite{huber2006b} still holds.

Now to the random bits.  Drawing a uniform number from $\{1,\ldots,n\}$
takes $\lceil \log_2 n \rceil$ bits, while drawing from $\{0,1\}$ for
coin $C_1$ takes one bit.  The expected number of bits needed to draw
a Bernoulli random variable with parameter not equal to $1/2$ 
is two, and so the total bits needed for
one step of the process (in expectation) is $\lceil \log_2 n \rceil + 3$.
Each step in the bounding chain and state uses at most two comparisons.
\end{proof}

It will be helpful in proving Theorem~\ref{THM:tpa} to have the following
bound on the tail of the Poisson distribution.

\begin{lemma}
\label{LEM:chernoff}
For $X \sim \pois(\mu)$ and $a \leq \mu$,
 $\prob(X \geq \mu + a) \leq \exp(-(1/2)a^2/\mu + (1/2)a^3/\mu^2)$ and 
for $a \leq \mu$, $\prob(X \leq \mu - a) \leq \exp(-a^2/(2\mu))$.
\end{lemma}

\begin{proof}
These follow from Chernoff Bounds which are essentially Markov's inequality
applied to the moment generating function of the random variable.  The
moment generating function of $X$ is $\mean[\exp(t X)] = \exp(\mu(e^t - 1)).$
So for $a > 0$
\[
\prob(X \geq \mu + a) = \prob(\exp(t X) \geq \exp(t(\mu + a)) 
 \leq \frac{\exp(\mu(e^t - 1))}{\exp(t(\mu + a))}.
\]
Setting $t = \ln(1 + a/\mu)$ minimizes the right hand side, and 
yields:
\[
\prob(X \geq \mu + a) \leq \exp(a - (\mu + a)\ln(1 + a/\mu)).
\]
For $a \leq \mu$, $-\ln(1 + a/\mu) \leq -a/\mu + (1/2)a^2/\mu^2$, 
so $\prob(X \geq \mu + a) \leq \exp(-(1/2)a^2/\mu + (1/2)a^3/\mu^2)$
as desired. For the second result:
\[
\prob(X \leq \mu - a) = \prob(\exp(-t X) \geq \exp(-t(\mu - a)) 
 \leq \frac{\exp(\mu(e^{-t} - 1))}{\exp(-t(\mu - a))}.
\]
Setting $t = -\ln(1 - a/\mu)$ then yields the next result. 
 
For $a \leq \mu,$ $-\ln(1 - a/\mu) \leq a/\mu + (1/2)(a/\mu)^2.$
So 
\[
-a - (\mu - a)\ln(1 - a/\mu) \leq -a + (\mu - a)((a/\mu) + (1/2)(a/\mu))
 = -(1/2)(a^2/\mu) - (1/2)(a^3/\mu^2).
\]
The right hand side is at most  $-(1/2)a^2/\mu,$ which completes the proof.
\end{proof}

\begin{proof}[Proof of Theorem~\ref{THM:tpa}]
Consider the first phase of the algorithm, where {\tpa} is run with 
$r_1 = 2\ln(2/\delta)$.  Consider the probability of the event 
$\{\hat A_1 + \sqrt{\hat A_1} + 2 < A\}.$
This event cannot happen if $A \leq 2$.  If $A > 2$, then this event occurs when
$\hat A_1 < A - (3/2) - \sqrt{A - 7/4}$.  Since 
$r_1 \hat A_1 \sim \pois(r_1A)$, 
Lemma~\ref{LEM:chernoff} can be used to say that
\begin{align*}
\prob(r_1 \hat A_1 < r_1A - r_1(3/2 + \sqrt{A - 7/4}))
 &\leq \exp(-(1/2)(r_1(3/2 + \sqrt{A - 7/4}))^2/(r_1 A) \\
 &\leq \exp(-(1/2)r_1 (9/4 + 3\sqrt{A - 7/4} + A - 7/4) / A) \\
 &\leq \exp(-(1/2)r_1) \\
 &\leq 2/\delta.
\end{align*}
In other words, with probability at least $1 - \delta/2$, 
$\hat A_1 + \sqrt{\hat A_1} + 2 \geq A$.

Now consider the second phase.  To simplify the notation, let
$\epsilon' = \ln(1 + \epsilon),$ and $\hat A_2 = \exp(\hat L_2)$ where
$\hat L_2$ is the output from the second phase.   Then from the first phase
$r_2 \geq A(\epsilon'^2 - \epsilon'^3)^{-1}\ln(4/\delta)$ with probability at
least $1 - \delta/2$.

So from Lemma~\ref{LEM:chernoff},
\begin{align*}
\prob(r_2 \hat A_2 \geq r_2 A + r_2 \epsilon') 
 &\leq \exp(-(1/2)(r_2 \epsilon')^2 / (r_2 A) 
       + (1/2) (r_2 \epsilon')^3/(r_2 A)^2) \\
 &= \exp(-(1/2)r_2 \epsilon'^2 / A + (1/2) r_2 \epsilon'^3/A^2) \\
 &\leq \exp(-\ln(4/\delta)).
\end{align*}
A similar bound holds for the left tail:
\[
\prob(r_2 \hat A_2 \leq r_2 A - r_2 \epsilon') 
 \leq \exp(-(1/2)r_2^2 \epsilon'^2/(r_2 A))
 \leq \delta / 4.
\]
Therefore, the total probability that failure occurs in either the
first phase or the second is at most 
$\delta / 2 + \delta / 4 + \delta / 4 = \delta.$  If $r_2 \hat A_2$
is within additive error $r_2 \epsilon' = r_2 \ln(1 + \epsilon)$ of $r_2 A$, 
then 
$\hat L_2 = \exp(\hat A_2/r)$ is within a factor of $1 + \epsilon$ of $\exp(A)$,
showing the result.
\end{proof}

To bound the expected running time, the following loose bound on the expected
value of the square root of a Poisson random variable is useful.
\begin{lemma}
\label{LEM:sqrt}
For $X \sim \pois(\mu)$, $\mean[\sqrt{X}] \leq \sqrt{\mu}.$
\end{lemma}

\begin{proof}
  Since $\sqrt{x}$ is a concave function, this follows from
  Jensen's inequality.
\end{proof}

\begin{proof}[Proof of Theorem~\ref{THM:bits}]
From Theorem~\ref{THM:nmcftp}, the expected number of bits per sample
is bounded by $4.3n^3(\ln n)(\lceil \log_2 n\rceil + 3)$ and does not
depend on the sample.  Hence the total number of expected bits can be
bounded by the expected number of bits per samples times the expected
number of samples.  The first phase of {\tpa} 
uses $r_1 = 2 \ln (2 / \delta)$ runs, each
with an expectation of $A + 1$ samples per run to make 
$r_1 (A + 1)$ expected samples.
The second phase uses 
$r_2 = (\hat A_1 + \sqrt{\hat A_1} + 2) 
 [\ln(1 + \epsilon)^2 - \ln(1 + \epsilon)^3] \ln(4/\delta)$ 
runs, where 
$r_1 \hat A_1 \sim \pois(r_1 A).$  So from 
Lemma~\ref{LEM:sqrt}, 
\[
\mean[\sqrt{\hat A_1}] = r_1^{-1/2} \mean[\sqrt{r_1 A}]
 \leq r_1^{-1/2} \sqrt{r_1 A} = \sqrt{A}.
\]
Using $A = \ln(\#(\linext))$ and then combining these factors yields the result.
\end{proof}

\section{Conclusion}
\label{SEC:conclusion}

{\tpa} is a sharp improvement on the self-reducibility method of 
Jerrum et al. for estimating the size of a set.  At first glance, the 
continuity requirement of {\tpa} precludes its use for discrete problems
such as linear extensions.  Fortunately, discrete problems can usually
be embedded in a continuous space to make the use of {\tpa} possible.
Here we have shown how to accomplish this task in such a way that the
time needed to take samples is the same as for uniform generation.
The result is 
an algorithm that is much faster
at estimating the number of linear extensions than previously known
algorithms.


\appendix

\section{R code}

The following R code implements these algorithms.

\begin{verbatim}
chain.step <- function(state,beta,i,c1,c2,posetmatrix) {
  # Takes one step in the Karzanov-Khatchyan chain using randomness in i, c1, c2
  # Assumes that home state is the identity permutation

  # Line 2 of Algorithm 5.1
  d <- i + 1 - state[i]
  # Line 3 of Algorithm 5.1
  if ((state[i] <= length(state)) && (state[i+1] <= length(state)))
    posetflag <- posetmatrix[state[i],state[i+1]]
  else
    posetflag <- 0
  if ( (c1 == 1) && (posetflag == 0) && 
       (d <= ceiling(beta)) && ((d < beta) || (c2 == 1))) {
      a <- state[i+1]; state[i+1] <- state[i]; state[i] <- a
    }

  return(state)
}

bounding.chain.step <- function(cstate,beta,i,c1,c2,posetmatrix) {
  # Based on Algorithm 5.2
  # Here cstate is a matrix with two rows, the first is the underlying state,
  # while the second is the bounding state
  
  n <- dim(cstate)[2]
  # Line 1
  if (cstate[2,i+1] == cstate[1,i]) c3 <- 1 - c1
  else c3 <- c1
  # Line 2 & 3
  cstate[1,] <- chain.step(cstate[1,],beta,i,c3,c2,posetmatrix)
  cstate[2,] <- chain.step(cstate[2,],beta,i,c1,c2,posetmatrix)
  # Line 4 through 6
  if (cstate[2,n] == (n+1)) 
    cstate[2,n] <- 1 + sum(cstate[2,] <= n)
  #Line 7
  return(cstate)
}

generate <- function(t,beta,posetmatrix) {
  n <- dim(posetmatrix)[1]; beta <- min(beta,n-1)
  # Line 1
  sigma <- 1:n; B <- rep(n+1,n)
  # Line 2 thorugh 4
  for (item in 1:(n - ceiling(beta)))
    B[item+ceiling(beta)] <- item
  B0 <- B
  cstate <- matrix(c(sigma,B),byrow=TRUE,nrow=2)
  # Line 5 through 8
  i <- floor(runif(t)*(n-1))+1;c1 <- rbinom(t,1,1/2)
  c2 <- rbinom(t,1,1+beta-ceiling(beta))
  for (s in 1:t) {
    cstate <- bounding.chain.step(cstate,beta,i[s],c1[s],c2[s],posetmatrix)
  }
  # Line 9
  if (sum(cstate[2,] < (n+1)) == n) return(cstate[2,])
  # Line 11-15
  else {
    cstate[1,] <- generate(2*t,beta,posetmatrix)
    cstate[2,] <- B0
    for (s in 1:t)
      cstate <- bounding.chain.step(cstate,beta,i[s],c1[s],c2[s],posetmatrix)
  }
  return(cstate[1,])
}

approximate.sample <- function(n,beta,posetmatrix,tvdist) {
  # Generates an approximate sample from the target distribution

  x <- 1:n
  n <- length(x)
  t <- 10*n^3*log(n)*log(1/tvdist)
  for (i in 1:t) {
    i <- runif(1)*(n-1)+1
    c1 <- rbinom(1,1,1/2)
    c2 <- rbinom(1,1,1+beta-ceiling(beta))
    x <- chain.step(x,beta,i,c1,c2,posetmatrix)
  }

  return(x)

}

checksum <- function(x) {
  
  checksum <- 0
  n <- length(x)
  for (i in n:1) {
    onespot <- which(x == 1)
    checksum <- checksum + factorial(i-1)*(onespot-1)
    x <- x[-onespot] - 1
  }
  return(checksum+1)
}

count.perfect.linear.extensions <- function(n = 4,beta = 4,posetmatrix,trials = 100) {
  # Generates a number of linear extensions, then counts the results
  
  results <- rep(0,factorial(n))
  # Burnin to an approximate sample
  x <- 1:n
  n <- length(x)
  # Take data
  for (i in 1:trials) {
    x <- generate(1,beta,posetmatrix)
    cs <- checksum(x)
    results[cs] <- results[cs] + 1
  }
  return(results/trials)
  
}

tpa.count.linear.extensions <- function(r,posetmatrix) {
  # Algorithm 6.1
  # Returns an estimate of the number of linear
  # extensions consistant with posetmatrix
  
  require(Matrix)
  
  n <- dim(posetmatrix)[1]
  # Line 1
  k <- 0
  # Line 2 through 12
  for (i in 1:r) {
    beta <- n - 1; k <- k - 1
    while (beta > 0) {
      k <- k + 1
      x <- generate(1,beta,posetmatrix) 
      xinv <- invPerm(x)
      betastep <- rep(0,n); y <- rep(0,n)
      for (j in 1:n) {
        y[j] <- runif(1)*((1+beta-ceiling(beta))*(xinv[j]-j == ceiling(beta))+
                        ((xinv[j]-j) < ceiling(beta)))
        betastep[j] <- (xinv[j]-j-1+y[j])*(xinv[j]-j > 0)
      }
      beta <- max(betastep)
 #     cat(" X: ",x,"\n X^{-1}: ",xinv,"\n Y: ",y,"\n betastep: ",betastep,"\n beta: ",beta,"\n")
    }
  }
  cat(" Estimate: [",exp((k-2*sqrt(k))/r),",",exp((k+2*sqrt(k))/r),"]\n")
  
  return(exp(k/r))
}

tpa.approximation <- function(posetmatrix,epsilon,delta) {
  # Gives an $(\epsilon,\delta)$-ras for the number of posets
  
  r1 <- ceiling(2*log(2/delta))
  a1 <- tpa.count.linear.extensions(r1,posetmatrix)
  a1 <- log(a1)
  r2 <- ceiling(2*(a1+sqrt(a1)+2)*
                  (log(1+epsilon)^2-log(1+epsilon)^3)^(-1)*log(4/delta))
  a2 <- tpa.count.linear.extensions(r2,posetmatrix)
  
  return(a2)
  
}

brute.force.count.linear.extensions <- function(posets) {
  # Counts the number of linear extensions of a poset by direct
  # ennumeration of all n! permutations and checking each to see
  # if it is a linear extension
  #
  # The poset is given as an n by n matrix whose (i,j)th entry
  # is the indicator function of $i \preceq j$
  
  require(gtools)
  
  n <- dim(posets)[1]
  A <- permutations(n,n)
  nfact <- nrow(A)
  le.flag <- rep(1,nfact)
  count <- 0
  for (i in 1:nfact) {
    for (a in 1:(n-1))
      for (b in (a+1):n) {
        le.flag[i] <- le.flag[i]*(1-posets[A[i,b],A[i,a]])
      }
    count <- count + le.flag[i]
  }
  return(count)
}

poset1 <- matrix(c(1,0,1,0,0,1,0,1,0,0,1,0,0,0,0,1),byrow=TRUE,ncol=4)
poset2 <- matrix(c(1,0,1,1,1,1,1,1, 0,1,0,1,0,1,1,1, 0,0,1,0,1,1,0,1,
                    0,0,0,1,0,1,1,1, 0,0,0,0,1,0,0,0, 0,0,0,0,0,1,0,1,
                    0,0,0,0,0,0,1,1, 0,0,0,0,0,0,0,1),byrow=TRUE,ncol=8)
poset3 <- t(matrix(c(1,1,0,1,0,1,0,1,0,0,1,1,0,0,0,1),nrow=4))
poset4 <- t(matrix(c(1,0,1,1,1,1, 0,1,0,1,1,1, 0,0,1,0,1,0,
                     0,0,0,1,1,1, 0,0,0,0,1,0, 0,0,0,0,0,1),nrow=6))
\end{verbatim}

\end{document}